\documentclass[11pt]{amsart}
\usepackage{amssymb,mathrsfs}
\title{A Partial Cayley Transform of Siegel-Jacobi Disk}

\author{Jae-Hyun Yang}
\address{Department of Mathematics, Inha University,
Incheon 402-751, Korea}
\email{jhyang@inha.ac.kr }

\begin{document}

\newtheorem{theorem}{Theorem}[section]
\newtheorem{lemma}{Lemma}[section]
\newtheorem{proposition}{Proposition}[section]
\newtheorem{remark}{Remark}[section]
\newtheorem{definition}{Definition}[section]

\renewcommand{\theequation}{\thesection.\arabic{equation}}
\renewcommand{\thetheorem}{\thesection.\arabic{theorem}}
\renewcommand{\thelemma}{\thesection.\arabic{lemma}}
\newcommand{\BR}{\mathbb R}
\newcommand{\BQ}{\mathbb Q}
\newcommand{\bn}{\bf n}
\def\charf {\mbox{{\text 1}\kern-.24em {\text l}}}
\newcommand{\BC}{\mathbb C}
\newcommand{\BZ}{\mathbb Z}


\begin{abstract}
Let ${\mathbb H}_g$ and ${\mathbb D}_g$ be the Siegel upper half
plane and the generalized unit disk of degree $g$ respectively.
Let ${\mathbb C}^{(h,g)}$ be the Euclidean space of all $h\times
g$ complex matrices. We present a partial Cayley transform of the
Siegel-Jacobi disk ${\mathbb D}_g\times {\mathbb C}^{(h,g)}$ onto
the Siegel-Jacobi space ${\mathbb H}_g\times {\mathbb C}^{(h,g)}$
which gives a partial bounded realization of ${\mathbb H}_g\times
{\mathbb C}^{(h,g)}$ by ${\mathbb D}_g\times {\mathbb C}^{(h,g)}$.
We prove that the natural actions of the Jacobi group on ${\mathbb
D}_g\times {\mathbb C}^{(h,g)}$ and ${\mathbb H}_g\times {\mathbb
C}^{(h,g)}$ are compatible via a partial Cayley transform. A
partial Cayley transform plays an important role in computing
differential operators on the Siegel-Jacobi disk ${\mathbb
D}_g\times {\mathbb C}^{(h,g)}$ invariant under the natural action
of the Jacobi group on ${\mathbb D}_g\times {\mathbb C}^{(h,g)}$
explicitly.
\end{abstract}
\maketitle

\newcommand\tr{\triangleright}
\newcommand\al{\alpha}
\newcommand\be{\beta}
\newcommand\g{\gamma}
\newcommand\gh{\Cal G^J}
\newcommand\G{\Gamma}
\newcommand\de{\delta}
\newcommand\e{\epsilon}
\newcommand\z{\zeta}
\newcommand\vth{\vartheta}
\newcommand\vp{\varphi}
\newcommand\om{\omega}
\newcommand\p{\pi}
\newcommand\la{\lambda}
\newcommand\lb{\lbrace}
\newcommand\lk{\lbrack}
\newcommand\rb{\rbrace}
\newcommand\rk{\rbrack}
\newcommand\s{\sigma}
\newcommand\w{\wedge}
\newcommand\fgj{{\frak g}^J}
\newcommand\lrt{\longrightarrow}
\newcommand\lmt{\longmapsto}
\newcommand\lmk{(\lambda,\mu,\kappa)}
\newcommand\Om{\Omega}
\newcommand\ka{\kappa}
\newcommand\ba{\backslash}
\newcommand\ph{\phi}
\newcommand\M{{\Cal M}}
\newcommand\bA{\bold A}
\newcommand\bH{\bold H}

\newcommand\Hom{\text{Hom}}
\newcommand\cP{\Cal P}
\newcommand\cH{\Cal H}

\newcommand\pa{\partial}

\newcommand\pis{\pi i \sigma}
\newcommand\sd{\,\,{\vartriangleright}\kern -1.0ex{<}\,}
\newcommand\wt{\widetilde}
\newcommand\fg{\frak g}
\newcommand\fk{\frak k}
\newcommand\fp{\frak p}
\newcommand\fs{\frak s}
\newcommand\fh{\frak h}
\newcommand\Cal{\mathcal}

\newcommand\fn{{\frak n}}
\newcommand\fa{{\frak a}}
\newcommand\fm{{\frak m}}
\newcommand\fq{{\frak q}}
\newcommand\CP{{\mathcal P}_g}
\newcommand\Hgh{{\mathbb H}_g \times {\mathbb C}^{(h,g)}}
\newcommand\Dgh{{\mathbb D}_g \times {\mathbb C}^{(h,g)}}
\newcommand\BD{\mathbb D}
\newcommand\BH{\mathbb H}
\newcommand\CCF{{\mathcal F}_g}
\newcommand\CM{{\mathcal M}}
\newcommand\Ggh{\Gamma_{g,h}}
\newcommand\Chg{{\mathbb C}^{(h,g)}}
\newcommand\Yd{{{\partial}\over {\partial Y}}}
\newcommand\Vd{{{\partial}\over {\partial V}}}

\newcommand\Ys{Y^{\ast}}
\newcommand\Vs{V^{\ast}}
\newcommand\LO{L_{\Omega}}
\newcommand\fac{{\frak a}_{\mathbb C}^{\ast}}
\newcommand\OW{\overline{W}}
\newcommand\OP{\overline{P}}
\newcommand\OQ{\overline{Q}}
\newcommand\Dg{{\mathbb D}_g}
\newcommand\Hg{{\mathbb H}_g}

%
%
\begin{section}{{\bf Introduction}}
\setcounter{equation}{0} For a given fixed positive integer $g$,
we let
$${\mathbb H}_g=\,\Big\{\,\Omega\in \BC^{(g,g)}\,\big|\ \Om=\,^t\Om,\ \ \ \text{Im}\,\Om>0\,\Big\}$$
be the Siegel upper half plane of degree $g$ and let
$$Sp(g,\BR)=\Big\{ M\in \BR^{(2g,2g)}\ \big\vert \ ^t\!MJ_gM= J_g\ \Big\}$$
be the symplectic group of degree $g$, where $F^{(k,l)}$ denotes
the set of all $k\times l$ matrices with entries in a commutative
ring $F$ for two positive integers $k$ and $l$, $^t\!M$ denotes
the transpose matrix of a matrix $M$ and
 $$J_g=\begin{pmatrix} 0&I_g\\
                   -I_g&0\end{pmatrix}.$$
We see that $Sp(g,\BR)$ acts on $\BH_g$ transitively by
\begin{equation}
M\cdot\Om=(A\Om+B)(C\Om+D)^{-1},
\end{equation} where $M=\begin{pmatrix} A&B\\
C&D\end{pmatrix}\in Sp(g,\BR)$ and $\Om\in \BH_g.$

\vskip 0.1cm Let $$\BD_g=\left\{\,W\in\BC^{(g,g)}\,|\ W=\,{ }^tW,\
I_g-W{\overline W}>0\,\right\}$$ be the generalized unit disk of
degree $g$. The Cayley transform $\Phi:\Dg\lrt\Hg$ defined by
\begin{equation}
\Phi(W)=i\,(I_g+W)(I_g-W)^{-1},\quad W\in\Dg
\end{equation}
is a biholomorphic mapping of $\Dg$ onto $\Hg$ which gives the
bounded realization of $\Hg$ by $\Dg$\,(cf.
\cite[pp.\,281-283]{Sie}). And the action (2.8) of the symplectic
group on $\BD_g$ is compatible with the action (1.1) via the
Cayley transform $\Phi$.

 \vskip 0.1cm For two positive integers $g$ and $h$, we consider
the Heisenberg group
$$H_{\BR}^{(g,h)}=\Big\{\,(\la,\mu;\ka)\,\big|\ \la,\mu\in \BR^{(h,g)},\ \kappa\in\BR^{(h,h)},\ \
\ka+\mu\,^t\la\ \text{symmetric}\ \Big\}$$ endowed with the
following multiplication law
$$(\la,\mu;\ka)\circ (\la',\mu';\ka')=(\la+\la',\mu+\mu';\ka+\ka'+\la\,^t\mu'-
\mu\,^t\la').$$ The Jacobi group $G^J$ is defined as the
semidirect product of $Sp(g,\BR)$ and $H_{\BR}^{(g,h)}$
$$G^J=Sp(g,\BR)\ltimes H_{\BR}^{(g,h)}$$
endowed with the following multiplication law
$$
\Big(M,(\lambda,\mu;\kappa)\Big)\cdot\Big(M',(\lambda',\mu';\kappa')\Big)
=\, \Big(MM',(\tilde{\lambda}+\lambda',\tilde{\mu}+ \mu';
\kappa+\kappa'+\tilde{\lambda}\,^t\!\mu'
-\tilde{\mu}\,^t\!\lambda')\Big)$$ with $M,M'\in Sp(g,\BR),
(\lambda,\mu;\kappa),\,(\lambda',\mu';\kappa') \in
H_{\BR}^{(g,h)}$ and
$(\tilde{\lambda},\tilde{\mu})=(\lambda,\mu)M'$. Then $G^J$ acts
on $\BH_g\times \BC^{(h,g)}$ transitively by
\begin{equation}
\Big(M,(\lambda,\mu;\kappa)\Big)\cdot
(\Om,Z)=\Big(M\cdot\Om,(Z+\lambda \Om+\mu)
(C\Omega+D)^{-1}\,\Big), \end{equation} where $M=\begin{pmatrix} A&B\\
C&D\end{pmatrix} \in Sp(g,\BR),\ (\lambda,\mu; \kappa)\in
H_{\BR}^{(g,h)}$ and $(\Om,Z)\in \BH_g\times \BC^{(h,g)}.$ In
\cite[p.\,1331]{YJH2}, the author presented the natural
construction of the action (1.3).

\vskip 0.2cm We mention that studying the Siegel-Jacobi space or
the Siegel-Jacobi disk associated with the Jacobi group is useful
to the study of the universal family of polarized abelian
varieties\,(cf.\,[10],\,[12]). The aim of this paper is to present
a partial Cayley transform of the Siegel-Jacobi disk $\Dgh$ onto
the Siegel-Jacobi space $\Hgh$ which gives a partially bounded
realization of $\Hgh$ by $\Dgh$ and to prove that the natural
actions of the Jacobi group on ${\mathbb D}_g\times {\mathbb
C}^{(h,g)}$ and ${\mathbb H}_g\times {\mathbb C}^{(h,g)}$ are
compatible via a partial Cayley transform. The main reason that we
study a partial Cayley transform is that this transform is
usefully applied to computing differential operators on the
Siegel-Jacobi disk $\Dgh$ invariant under the action (3.5) of the
Jacobi group $G_*^J$\,(cf.\,(3.2)) explicitly.

\vskip 0.2cm This paper is organized as follows. In Section 2, we
review the Cayley transform of the generalized unit disk $\Dg$
onto the Siegel upper half plane $\BH_g$ which gives a bounded
realization of ${\mathbb H}_g$ by ${\mathbb D}_g$ . In Section 3,
we construct a partial Cayley transform of the Siegel-Jacobi disk
$\Dgh$ onto the Siegel-Jacobi space $\Hgh$ which gives a partially
bounded realization of $\Hgh$ by $\Dgh$\,(cf.\,(3.6)). We prove
that the action (1.3) of the Jacobi group $G^J$ is compatible with
the action (3.5) of the Jacobi group $G_*^J$ through a partial
Cayley transform (cf. Theorem 3.1).  In the final section, we
present the canonical automorphic factors of the Jacobi group
$G^J_*$.

\vskip 0.2cm  \noindent $ \textsc{Notations\,:}$  \ We denote by
$\BR$ and $\BC$ the field of real numbers, and the field of
complex numbers respectively. For a square matrix $A\in F^{(k,k)}$
of degree $k$, $\sigma(A)$ denotes the trace of $A$. For
$\Omega\in {\mathbb H}_g,\ \text{Re}\,\Omega$ ({\it resp.}\
$\textrm{Im}\,\Omega)$ denotes the real ({\it resp.}\ imaginary)
part of $\Omega.$  For a matrix $A\in F^{(k,k)}$ and $B\in
F^{(k,l)},$ we write $A[B]=\,^tBAB$. $I_n$ denotes the identity
matrix of degree $n$.

\end{section}

%
%
\begin{section}{{\bf The Cayley Transform }}
\setcounter{equation}{0}

\vskip 0.2cm Let
\begin{equation}
T={1\over {\sqrt{2}} }\,
\begin{pmatrix} \ I_g&\ I_g\\
                   iI_g&-iI_g\end{pmatrix}
\end{equation}
be the $2g\times 2g$ matrix represented by $\Phi.$ Then
\begin{equation}
T^{-1}Sp(g,\BR)\,T=\left\{ \begin{pmatrix} P & Q \\ \OQ & \OP
\end{pmatrix}\,\Big|\ ^tP\OP-\,{}^t\OQ Q=I_g,\ {}^tP\OQ=\,{}^t\OQ
P\,\right\}.
\end{equation}
Indeed, if $M=\begin{pmatrix} A&B\\
C&D\end{pmatrix}\in Sp(g,\BR)$, then
\begin{equation}
T^{-1}MT=\begin{pmatrix} P & Q \\ \OQ & \OP
\end{pmatrix},
\end{equation}
where
\begin{equation}
P= {\frac 12}\,\Big\{ (A+D)+\,i\,(B-C)\Big\}
\end{equation}
and
\begin{equation}
 Q={\frac
12}\,\Big\{ (A-D)-\,i\,(B+C)\Big\}.
\end{equation}

For brevity, we set
\begin{equation*}
G_*= T^{-1}Sp(g,\BR)T.
\end{equation*}
Then $G_*$ is a subgroup of $SU(g,g),$ where
$$SU(g,g)=\left\{\,h\in\BC^{(g,g)}\,\big|\ {}^th I_{g,g}{\overline
h}=I_{g,g}\,\right\},\quad I_{g,g}=\begin{pmatrix} \ I_g&\ 0\\
0&-I_g\end{pmatrix}.$$ In the case $g=1$, we observe that
$$T^{-1}Sp(1,\BR)T=T^{-1}SL_2(\BR)T=SU(1,1).$$
If $g>1,$ then $G_*$ is a {\it proper} subgroup of $SU(g,g).$ In
fact, since ${}^tTJ_gT=-\,i\,J_g$, we get
\begin{equation}G_*=\Big\{\,h\in SU(g,g)\,\big|\
{}^thJ_gh=J_g\,\Big\}=SU(g,g)\cap Sp(g,\BC),
\end{equation}

\noindent where
$$Sp(g,\BC)=\Big\{\alpha\in \BC^{(2g,2g)}\ \big\vert \ ^t\!\alpha\, J_g\,\alpha= J_g\ \Big\}.$$

Let
\begin{equation*}
P^+=\left\{\begin{pmatrix} I_g & Z\\ 0 & I_g
\end{pmatrix}\,\Big|\ Z=\,{}^tZ\in\BC^{(g,g)}\,\right\}
\end{equation*}
be the $P^+$-part of the complexification of $G_*\subset SU(g,g).$
We note that the Harish-Chandra decomposition of an element
$\left(\begin{pmatrix} P & Q\\ {\overline Q} & {\overline P}
\end{pmatrix}\right)$ in $G_*$ is
\begin{equation*}
\begin{pmatrix} P & Q\\ \OQ & \OP
\end{pmatrix}=\begin{pmatrix} I_g & Q\OP^{-1}\\ 0 & I_g
\end{pmatrix} \begin{pmatrix} P-Q\OP^{-1}\OQ & 0\\ 0 & \OP
\end{pmatrix} \begin{pmatrix} I_g & 0\\ \OP^{-1}\OQ & I_g
\end{pmatrix}.
\end{equation*}
For more detail, we refer to [2,\,\,p.\,155]. Thus the
$P^+$-component of the following element
$$\begin{pmatrix} P & Q\\ \OQ & \OP
\end{pmatrix}   \cdot\begin{pmatrix} I_g & W\\ 0 & I_g
\end{pmatrix},\quad W\in \BD_g$$ of the complexification of $G_*^J$ is
given by
\begin{equation}
\left( \begin{pmatrix} I_g & (PW+Q)(\OQ W+\OP)^{-1}
\\ 0 & I_g
\end{pmatrix}\right).
\end{equation}
\newcommand\POB{ {{\partial}\over {\partial{\overline \Omega}}} }
\newcommand\PZB{ {{\partial}\over {\partial{\overline Z}}} }
\newcommand\PX{ {{\partial}\over{\partial X}} }
\newcommand\PY{ {{\partial}\over {\partial Y}} }
\newcommand\PU{ {{\partial}\over{\partial U}} }
\newcommand\PV{ {{\partial}\over{\partial V}} }
\newcommand\PO{ {{\partial}\over{\partial \Omega}} }
\newcommand\PZ{ {{\partial}\over{\partial Z}} }
\newcommand\PW{ {{\partial}\over{\partial W}} }
\newcommand\PWB{ {{\partial}\over {\partial{\overline W}}} }
\newcommand\OVW{\overline W}

\noindent We note that $Q\OP^{-1}\in\Dg.$ We get the
Harish-Chandra embedding of $\Dg$ into $P^+$\,(cf.
\cite[p.\,155]{Kn} or \cite[pp.\,58-59]{Sat}). Therefore we see
that $G_*$ acts on $\Dg$ transitively by
\begin{equation}
\begin{pmatrix} P & Q \\ \OQ & \OP
\end{pmatrix}\cdot W=(PW+Q)(\OQ W+\OP)^{-1},\quad \begin{pmatrix} P & Q \\ \OQ & \OP
\end{pmatrix}\in G_*,\ W\in \Dg.
\end{equation}
The isotropy subgroup at the origin $o$ is given by
$$K=\left\{\,\begin{pmatrix} P & 0 \\ 0 & {\overline
P}\end{pmatrix}\,\Big|\ P\in U(g)\ \right\}.$$ Thus $G_*/K$ is
biholomorphic to $\Dg$. It is known that the action (1.1) is
compatible with the action (2.8) via the Cayley transform $\Phi$\
(cf.\,(1.2)). In other words, if $M\in Sp(g,\BR)$ and $W\in\BD_g$,
then
\begin{equation}
M\cdot \Phi(W)=\Phi(M_*\cdot W),
\end{equation}

\noindent where $M_*=T^{-1}MT\in G_*.$ For a proof of Formula
(2.9), we refer to the proof of Theorem 3.1.

\vskip 0.2cm
 For $\Om=(\omega_{ij})\in\BH_g,$ we write $\Om=X+iY$
with $X=(x_{ij}),\ Y=(y_{ij})$ real and $d\Om=(d\om_{ij})$. We
also put
$$\PO=\,\left(\,
{ {1+\delta_{ij}}\over 2}\, { {\partial}\over {\partial \om_{ij} }
} \,\right) \qquad\text{and}\qquad \POB=\,\left(\, {
{1+\delta_{ij}}\over 2}\, { {\partial}\over {\partial {\overline
{\om}}_{ij} } } \,\right).$$ Then
\begin{equation}
ds^2=\s \Big( Y^{-1}d\Om\,
Y^{-1}d{\overline\Om}\Big)\end{equation} is a
$Sp(g,\BR)$-invariant metric on $\BH_g$ (cf.\,\cite{Sie}) and H.
Maass \cite{Mas} proved that its Laplacian is given by
\begin{equation}
\Delta=\,4\,\s \left( Y\,\,{
}^t\!\left(Y\POB\right)\PO\right).\end{equation}

\vskip 0.2cm For $W=(w_{ij})\in \Dg,$ we write $dW=(dw_{ij})$ and
$d{\overline W}=(d{\overline{w}}_{ij})$. We put $$\PW=\,\left(\, {
{1+\delta_{ij}}\over 2}\, { {\partial}\over {\partial w_{ij} } }
\,\right) \qquad\text{and}\qquad \PWB=\,\left(\, {
{1+\delta_{ij}}\over 2}\, { {\partial}\over {\partial {\overline
{w}}_{ij} } } \,\right).$$

Using the Cayley transform $\Phi:\Dg\lrt \BH_g$, H. Maass proved
(cf.\,\,[3]) that
\begin{equation}
ds_*^2=4 \s \Big((I_g-W{\overline W})^{-1}dW\,(I_g-\OVW
W)^{-1}d\OVW\,\Big)\end{equation} is a $G_*$-invariant Riemannian
metric on $\BD_g$ and its Laplacian is given by
\begin{equation}
\Delta_*=\,\s \left( (I_g-W\OW)\,{ }^t\!\left( (I_g-W\OVW)
\PWB\right)\PW\right).\end{equation}

\end{section}
%
%
\begin{section}{{\bf A Partial Cayley Transform}}
\setcounter{equation}{0}

\vskip 0.15cm In this section, we present a partial Cayley
transform of $\Dgh$ onto $\Hgh$ and prove that the action (1.3) of
$G^J$ on $\Hgh$ is compatible with the {\it natural
action}\,(cf.\,(3.5)) of the Jacobi group $G_*^J$ on $\Dgh$ via a
partial Cayley transform.

\vskip 0.2cm From now on, for brevity we write
$\BH_{g,h}=\BH_g\times \BC^{(h,g)}.$ We can identify an element
$g=\Big(M,(\la,\mu;\kappa)\Big)$
of $G^J,\ M=\begin{pmatrix} A&B\\
C&D\end{pmatrix}\in Sp(g,\BR)$ with the element
\begin{equation*}
\begin{pmatrix} A & 0 & B & A\,^t\mu-B\,^t\la  \\ \la & I_h & \mu
& \kappa \\ C & 0 & D & C\,^t\mu-D\,^t\la \\ 0 & 0 & 0 & I_h
\end{pmatrix}
\end{equation*}
of $Sp(g+h,\BR).$ This subgroup plays an important role in
investigating the Fourier-Jacobi expansion of a Siegel modular
form for $Sp(g+h,\BR)$\,(cf.\,[4]) and studying the theory of
Jacobi forms\,(cf.\,[1],\,[7-9],\,[15]).

\vskip 0.3cm We set
\begin{equation*}
T_*={1\over {\sqrt 2}}\,
\begin{pmatrix} I_{g+h} & I_{g+h}\\ iI_{g+h} & -iI_{g+h}
\end{pmatrix}.
\end{equation*}
We now consider the group $G_*^J$ defined by
\begin{equation*}
G_*^J=T_*^{-1}G^JT_*.
\end{equation*}
If $g=(M,(\la,\mu;\kappa))\in G^J$ with $M=\begin{pmatrix} A&B\\
C&D\end{pmatrix}\in Sp(g,\BR)$, then $T_*^{-1}gT_*$ is given by
\begin{equation}
T_*^{-1}gT_*=
\begin{pmatrix} P_* & Q_*\\ {\overline Q}_* & {\overline P}_*
\end{pmatrix},
\end{equation}
where
\begin{equation*}
P_*=
\begin{pmatrix} P & {\frac 12} \Big\{ Q\,\,{}^t(\la+i\mu)-P\,\,{}^t(\la-i\mu)\Big\}\\
{\frac 12} (\la+i\mu) & I_h+i{\frac \kappa 2}
\end{pmatrix},
\end{equation*}

\begin{equation*}
Q_*=
\begin{pmatrix} Q & {\frac 12} \Big\{ P\,\,{}^t(\la-i\mu)-Q\,\,{}^t(\la+i\mu)\Big\}\\
{\frac 12} (\la-i\mu) & -i{\frac \kappa 2}
\end{pmatrix},
\end{equation*}
and $P,\,Q$ are given by Formulas (2.4) and (2.5). From now on, we
write
\begin{equation*}
\left(\begin{pmatrix} P & Q\\ {\overline Q} & {\overline P}
\end{pmatrix},\left( {\frac 12}(\la+i\mu),\,{\frac 12}(\la-i\mu);\,-i{\kappa\over 2}\right)\right)=
\begin{pmatrix} P_* & Q_*\\ {\overline Q}_* & {\overline P}_*
\end{pmatrix}.
\end{equation*}
In other words, we have the relation
\begin{equation*}
T_*^{-1}\left( \begin{pmatrix} A&B\\
C&D\end{pmatrix},(\la,\mu;\kappa)
\right)T_*=  \left(\begin{pmatrix} P & Q\\
{\overline Q} & {\overline P}
\end{pmatrix},\left(
{\frac 12}(\la+i\mu),\,{\frac 12}(\la-i\mu);\,-i{\kappa\over 2}
\right)\right).
\end{equation*}
Let
\begin{equation*}
H_{\BC}^{(g,h)}=\left\{ (\xi,\eta\,;\zeta)\,|\
\xi,\eta\in\BC^{(h,g)},\ \zeta\in\BC^{(h,h)},\
\zeta+\eta\,{}^t\xi\ \textrm{symmetric}\,\right\}
\end{equation*}
be the Heisenberg group endowed with the following multiplication
\begin{equation*}
(\xi,\eta\,;\zeta)\circ
(\xi',\eta';\zeta')=(\xi+\xi',\eta+\eta'\,;\zeta+\zeta'+
\xi\,{}^t\eta'-\eta\,{}^t\xi')).
\end{equation*}
We define the semidirect product
\begin{equation*}
SL(2g,\BC)\ltimes H_{\BC}^{(g,h)}
\end{equation*}
endowed with the following multiplication
\begin{eqnarray*}
& & \left( \begin{pmatrix} P & Q\\ R & S
\end{pmatrix}, (\xi,\eta\,;\zeta)\right)\cdot \left( \begin{pmatrix} P' & Q'\\
R' & S'
\end{pmatrix}, (\xi',\eta';\zeta')\right)\\
&=& \left( \begin{pmatrix} P & Q\\ R & S
\end{pmatrix}\,\begin{pmatrix} P' & Q'\\ R' & S'
\end{pmatrix},\,({\tilde \xi}+\xi',{\tilde
\eta}+\eta';\zeta+\zeta'+{\tilde \xi}\,{}^t\eta'-{\tilde
\eta}\,{}^t\xi')  \right),
\end{eqnarray*}
where ${\tilde\xi}=\xi P'+\eta R'$ and ${\tilde \eta}=\xi Q'+\eta
S'.$

\vskip 0.2cm If we identify $H_{\BR}^{(g,h)}$ with the subgroup
$$\left\{ (\xi,{\overline \xi};i\kappa)\,|\ \xi\in\BC^{(h,g)},\
\ka\in\BR^{(h,h)}\,\right\}$$ of $H_{\BC}^{(g,h)},$ we have the
following inclusion
$$G_*^J\subset SU(g,g)\ltimes H_{\BR}^{(g,h)}\subset SL(2g,\BC)\ltimes
H_{\BC}^{(g,h)}.$$ More precisely, if we recall $G_*=SU(g,g)\cap
Sp(g,\BC)$ (cf.\,(2.6)), we see that the Jacobi group $G_*^J$ is
given by
\begin{equation}
G^J_*=\left\{ \left( \begin{pmatrix} P & Q\\ \OQ & \OP
\end{pmatrix}, (\xi,{\bar \xi}\,;i\kappa)\right)\,\Big|
\ \begin{pmatrix} P & Q\\ \OQ & \OP
\end{pmatrix}
\in G_*,\ \xi\in \BC^{(m,n)},\ \kappa\in \BR^{(m,m)}\,\right\}.
\end{equation}

\noindent We define the mapping $\Theta:G^J\lrt G_*^J$ by
\begin{equation}\Theta
\left( \begin{pmatrix} A&B\\
C&D\end{pmatrix},(\la,\mu;\kappa) \right)=\left(\begin{pmatrix} P
& Q\\ {\overline Q} & {\overline P}
\end{pmatrix},\left(
{\frac 12}(\la+i\mu),\,{\frac 12}(\la-i\mu)\,;\,-i{\kappa\over 2}
\right)\right),\end{equation} where $P$ and $Q$ are given by
Formulas (2.4) and (2.5). We can see that if $g_1,g_2\in G^J$,
then $\Theta(g_1g_2)=\Theta(g_1)\Theta(g_2).$

\vskip 0.2cm According to \cite[p.\,250]{YJH3}, $G_*^J$ is of the
Harish-Chandra type\,(cf. \cite[p.\,118]{Sat}). Let
$$g_*=\left(\begin{pmatrix} P & Q\\
{\overline Q} & {\overline P}
\end{pmatrix},\left( \la, \mu;\,\kappa\right)\right)$$
be an element of $G_*^J.$ Since the Harish-Chandra decomposition
of an element $\begin{pmatrix} P & Q\\ R & S
\end{pmatrix}$ in $SU(g,g)$ is given by
\begin{equation*}
\begin{pmatrix} P & Q\\ R & S
\end{pmatrix}=\begin{pmatrix} I_g & QS^{-1}\\ 0 & I_g
\end{pmatrix} \begin{pmatrix} P-QS^{-1}R & 0\\ 0 & S
\end{pmatrix} \begin{pmatrix} I_g & 0\\ S^{-1}R & I_g
\end{pmatrix},
\end{equation*}
the $P_*^+$-component of the following element
$$g_*\cdot\left( \begin{pmatrix} I_g & W\\ 0 & I_g
\end{pmatrix}, (0,\eta;0)\right),\quad W\in \BD_g$$
of $SL(2g,\BC)\ltimes H_{\BC}^{(g,h)}$ is given by
\begin{equation}
\left( \begin{pmatrix} I_g & (PW+Q)(\OQ W+\OP)^{-1}
\\ 0 & I_g
\end{pmatrix},\,\left(0,\,(\eta+\la W+\mu)(\OQ W+\OP)^{-1}\,;0\right)\right).
\end{equation}

\vskip 0.2cm We can identify $\Dgh$ with the subset
\begin{equation*}
\left\{ \left( \begin{pmatrix} I_g & W\\ 0 & I_g
\end{pmatrix}, (0,\eta;0)\right)\,\Big|\ W\in\BD_g,\
\eta\in\BC^{(h,g)}\,\right\}\end{equation*} of the
complexification of $G_*^J.$ Indeed, $\Dgh$ is embedded into
$P_*^+$ given by
\begin{equation*}
P_*^+=\left\{\,\left( \begin{pmatrix} I_g & W\\ 0 & I_g
\end{pmatrix}, (0,\eta;0)\right)\,\Big|\ W=\,{}^tW\in \BC^{(g,g)},\
\eta\in\BC^{(h,g)}\ \right\}.
\end{equation*}
This is a generalization of the Harish-Chandra embedding\,(cf.
\cite[p.\,119]{Sat}). Hence $G_*^J$ acts on $\Dgh$ transitively by
\begin{equation}
\left(\begin{pmatrix} P & Q\\
{\overline Q} & {\overline P}
\end{pmatrix},\left( \la, \mu;\,\kappa\right)\right)\cdot
(W,\eta)=\Big((PW+Q)(\OQ W+\OP)^{-1},(\eta+\la W+\mu)(\OQ
W+\OP)^{-1}\,\Big).\end{equation}

\noindent \vskip 0.2cm From now on, for brevity we write
$\BD_{g,h}=\Dgh.$ We define the map $\Phi_*$ of $\BD_{g,h}$ into
$\BH_{g,h}$ by
\begin{equation}
\Phi_*(W,\eta)=\Big(
i(I_g+W)(I_g-W)^{-1},\,2\,i\,\eta\,(I_g-W)^{-1}\Big),\quad
(W,\eta)\in\BD_{g,h}.
\end{equation}

\noindent We can show that $\Phi_*$ is a biholomorphic map of
$\BD_{g,h}$ onto $\BH_{g,h}$ which gives a partial bounded
realization of $\BH_{g,h}$ by the Siegel-Jacobi disk $\BD_{g,h}$.
We call this map $\Phi_*$ the $ \textit{partial Cayley transform}$
of the Siegel-Jacobi disk $\BD_{g,h}$.

\vskip 0.2cm

\begin{theorem}
The action (1.3) of $G^J$ on $\BH_{g,h}$ is compatible with the
action (3.5) of $G_*^J$ on $\BD_{g,h}$ through the {\it partial
Cayley transform} $\Phi_*$. In other words, if $g_0\in G^J$ and
$(W,\eta)\in\BD_{g,h}$,
\begin{equation}
g_0\cdot\Phi_*(W,\eta)=\Phi_*(g_*\cdot (W,\eta)),
\end{equation}
where $g_*=T_*^{-1}g_0 T_*$. We observe that Formula (3.7)
generalizes Formula (2.9). The inverse of $\Phi_*$ is
\begin{equation}
\Phi_*^{-1}(\Omega,Z)=\Big((\Omega-iI_g)(\Omega+iI_g)^{-1},\,Z(\Omega+iI_g)^{-1}\Big).
\end{equation}
\end{theorem}

\begin{proof}
 Let
 \begin{equation*}
 g_0=\left(\begin{pmatrix} A&B\\
C&D\end{pmatrix},(\la,\mu;\kappa) \right)
\end{equation*}
be an element of $G^J$ and let $g_*=T_*^{-1}g_0T_*.$ Then
 \begin{equation*}
g_*=\left(\begin{pmatrix} P & Q\\ {\overline Q} & {\overline P}
\end{pmatrix},\left( {\frac 12}(\la+i\mu),\,{\frac 12}(\la-i\mu)\,;\,-i{\kappa\over
2}\right)\right),
\end{equation*}
where $P$ and $Q$ are given by Formulas (2.4) and (2.5).

\vskip 0.2cm For brevity, we write
\begin{equation*}
(\Omega,Z)=\Phi_*(W,\eta)\quad \emph{and}\quad
(\Omega_*,Z_*)=g_0\cdot (\Omega,Z).
\end{equation*}
That is,
\begin{equation*}
\Omega=i(I_g+W)(I_g-W)^{-1}\quad \emph{and}\quad Z=2\,i\,\eta
(I_g-W)^{-1}.
\end{equation*}
Then we get
\begin{eqnarray*}
\Omega_*&=&(A\Omega+B)(C\Omega+D)^{-1}\\
&=&\Big\{ i\,A(I_g+W)(I_g-W)^{-1}+B\Big\}\Big\{
i\,C(I_g+W)(I_g-W)^{-1}+D\Big\}^{-1}\\
&=&\Big\{ i\,A(I_g+W)+B(I_g-W)\Big\}(I_g-W)^{-1}\\
& &\ \times \Big[   \Big\{
i\,C(I_g+W)+D(I_g-W)\Big\} (I_g-W)^{-1}\Big]^{-1}\\
&=&\Big\{ (i\,A-B)W+(i\,A+B)\Big\}\Big\{
(i\,C-D)W+(i\,C+D)\Big\}^{-1}
\end{eqnarray*}
and
\begin{eqnarray*}
Z_*&=&(Z+\la \Omega+\mu)(C\Omega+D)^{-1}\\
&=&\Big\{ 2\,i\, \eta\,(I_g-W)^{-1}+i\,\la
(I_g+W)(I_g-W)^{-1}+\mu\Big\}\\
& & \ \times \Big\{ i\,
C(I_g+W)(I_g-W)^{-1}+D\Big\}^{-1}\\
&=& \Big\{ 2\,i\,\eta+i\,\la\,(I_g+W)+\mu\, (I_g-W)\Big\}(I_g-W)^{-1}\\
& & \ \times \Big[
\Big\{ i\,C(I_g+W)+D(I_g-W)\Big\}(I_g-W)^{-1}\Big]^{-1}\\
&=& \Big\{ 2\,i\,\eta+(\la\, i-\mu)W+\la\, i+\mu\Big\}\Big\{
(i\,C-D)W+(i\,C+D)\Big\}^{-1}.
\end{eqnarray*}
On the other hand, we set
\begin{equation*}
(W_*,\eta_*)=g_*\cdot(W,\eta)\quad \emph{and}\quad
({\widehat\Omega},{\widehat Z})=\Phi_*(W_*,\eta_*).
\end{equation*}
Then
\begin{equation*}
W_*=(PW+Q)(\OQ W+\OP)^{-1}  \quad \emph{and}\quad
\eta_*=(\eta+\la_* W+\mu_*)(\OQ W+\OP)^{-1},
\end{equation*}
where $\la_*={\frac 12}(\la+\,i\,\mu)$ and $\mu_*={\frac
12}(\la-\,i\,\mu).$

\vskip 0.2cm According to Formulas (2.4) and (2.5), we get
\begin{eqnarray*}
{\widehat \Omega}&=&i\,(I_g+W_*)(I_g-W_*)^{-1}\\
&=&i\,\Big\{ I_g+(PW+Q)(\OQ W+\OP)^{-1}\Big\}\Big\{
I_g-(PW+Q)(\OQ W+\OP)^{-1}\Big\}^{-1}\\
&=& i\,(\OQ W+\OP+PW+Q)(\OQ W+\OP)^{-1}\\
& &\ \times \Big\{ (\OQ W+\OP-PW-Q)(\OQ W+\OP)^{-1}\Big\}^{-1}\\
&=& i\,\Big\{ (P+\OQ\,)W+\OP+Q\Big\}\Big\{
(\OQ-P)W+\OP-Q\Big\}^{-1}\\
&=&\Big\{ (i\,A-B)W+(i\,A+B)\Big\}\Big\{
(i\,C-D)W+(i\,C+D)\Big\}^{-1}
\end{eqnarray*}
Therefore ${\widehat \Omega}=\Omega_*.$ In fact, this result is
the well-known fact  (cf. Formula (2.9)) that the action (1.1) is
compatible with the action (2.8) via the Cayley transform $\Phi$.
Here we gave a proof of Formula (2.9) for convenience.
\begin{eqnarray*}
{\widehat Z}&=& 2\,i\,\eta_*(I_g-W_*)^{-1}\\
&=& 2\,i\,(\eta+\la_* W+\mu_* )(\OQ W+\OP)^{-1}\\
& & \ \times \Big\{ I_g-(PW+Q)(\OQ W+\OP\,)^{-1}\Big\}^{-1}\\
&=& 2\,i\,(\eta+\la_*W+\mu_*)(\OQ W+\OP\,)^{-1}\\
& &\ \times \Big\{ (\OQ W+\OP-PW-Q)(\OQ
W+\OP\,)^{-1}\Big\}^{-1}\\
&=& 2\,i\,(\eta+\la_* W+\mu_*)\Big\{ (\OQ-P)W+\OP-Q\Big\}^{-1}.
\end{eqnarray*}
Using Formulas (2.4) and (2.5), we obtain
\begin{equation*}
{\widehat Z}=\Big\{
2\,i\,\eta+(\la\,i-\mu)W+\la\,i+\mu\Big\}\Big\{
(i\,C-D)W+i\,C+D\Big\}^{-1}.
\end{equation*}
Hence ${\widehat Z}=Z_*.$ Consequently we get Formula (3.7).
Formula (3.8) follows immediately from a direct computation.
\end{proof}

\vskip 0.4cm For a coordinate $(\Om,Z)\in\BH_{g,h}$ with
$\Om=(\om_{\mu\nu})\in {\mathbb H}_g$ and $Z=(z_{kl})\in \Chg,$ we
put
\begin{align*}
\Om\,=&\,X\,+\,iY,\quad\ \ X\,=\,(x_{\mu\nu}),\quad\ \
Y\,=\,(y_{\mu\nu})
\ \ \text{real},\\
Z\,=&U\,+\,iV,\quad\ \ U\,=\,(u_{kl}),\quad\ \ V\,=\,(v_{kl})\ \
\text{real},\\
d\Om\,=&\,(d\om_{\mu\nu}),\quad\ \ dX\,=\,(dx_{\mu\nu}),\quad\ \
dY\,=\,(dy_{\mu\nu}),\\
dZ\,=&\,(dz_{kl}),\quad\ \ dU\,=\,(du_{kl}),\quad\ \
dV\,=\,(dv_{kl}),\\
d{\overline{\Om}}=&\,(d{\overline{\om}}_{\mu\nu}),\quad
d{\overline Z}=(d{\bar z}_{kl}),
\end{align*}
\newcommand\POB{ {{\partial}\over {\partial{\overline \Omega}}} }
\newcommand\PZB{ {{\partial}\over {\partial{\overline Z}}} }
\newcommand\PX{ {{\partial}\over{\partial X}} }
\newcommand\PY{ {{\partial}\over {\partial Y}} }
\newcommand\PU{ {{\partial}\over{\partial U}} }
\newcommand\PV{ {{\partial}\over{\partial V}} }
\newcommand\PO{ {{\partial}\over{\partial \Omega}} }
\newcommand\PZ{ {{\partial}\over{\partial Z}} }
\newcommand\bw{d{\overline W}}
\newcommand\bz{d{\overline Z}}
\newcommand\bo{d{\overline \Omega}}

$$ {{\partial}\over{\partial \Omega}}\,=\,\left(\, { {1+\delta_{\mu\nu}} \over 2}\, {
{\partial}\over {\partial \om_{\mu\nu}} } \,\right),\quad
\POB\,=\,\left(\, { {1+\delta_{\mu\nu}}\over 2} \, {
{\partial}\over {\partial {\overline \om}_{\mu\nu} }  }
\,\right),$$
$$\PZ=\begin{pmatrix} { {\partial}\over{\partial z_{11}} } & \hdots &
{ {\partial}\over{\partial z_{h1}} }\\
\vdots&\ddots&\vdots\\
{ {\partial}\over{\partial z_{1g}} }&\hdots &{ {\partial}\over
{\partial z_{hg}} } \end{pmatrix},\quad \PZB=\begin{pmatrix} {
{\partial}\over{\partial {\overline z}_{11} }   }&
\hdots&{ {\partial}\over{\partial {\overline z}_{h1} }  }\\
\vdots&\ddots&\vdots\\
{ {\partial}\over{\partial{\overline z}_{1g} }  }&\hdots & {
{\partial}\over{\partial{\overline z}_{hg} }  }
\end{pmatrix}.$$
\vskip 0.1cm\noindent $ \textbf{Remark 3.1.}$ The author proved in
\cite{YJH4} that for any two positive real numbers $A$ and $B$,
the following metric
\begin{eqnarray}
ds_{g,h;A,B}^2&=&\,A\, \s\Big( Y^{-1}d\Om\,Y^{-1}d{\overline
\Om}\Big) \nonumber \\ && \ \ + \,B\,\bigg\{ \s\Big(
Y^{-1}\,^tV\,V\,Y^{-1}d\Om\,Y^{-1} \bo \Big)
 +\,\s\Big( Y^{-1}\,^t(dZ)\,\bz\Big)\\
&&\quad\quad -\s\Big( V\,Y^{-1}d\Om\,Y^{-1}\,^t(\bz)\Big)\,
-\,\s\Big( V\,Y^{-1}\bo\, Y^{-1}\,^t(dZ)\,\Big) \bigg\} \nonumber
\end{eqnarray}

\noindent is a Riemannian metric on $\BH_{g,h}$ which is invariant
under the action (1.3) of the Jacobi group $G^J$ and its Laplacian
is given by
\begin{eqnarray}
\Delta_{n,m;A,B}\,&=& \frac4A \,\bigg\{ \s\left(\,Y\,\,
^t\!\left(Y\POB\right)\PO\,\right)\, +\,\s\left(\,VY^{-1}\,^tV\,\,^t\!\left(Y\PZB\right)\,\PZ\,\right)\nonumber\\
& &\ \
+\,\s\left(V\,\,^t\!\left(Y\POB\right)\PZ\,\right)+\,\s\left(\,^tV\,\,^t\!\left(Y\PZB\right)\PO\,\right)\bigg\}\\
& & \ +\frac4B\,\,\s\left(\, Y\,\PZ\,\,{}^t\!\left(
\PZB\right)\,\right).\nonumber
\end{eqnarray}
We observe that Formulas (3.9) and (3.10) generalize Formulas
(2.10) and (2.11). The following differential form
$$dv_{g,h}=\,\left(\,\text{det}\,Y\,\right)^{-(g+h+1)}[dX]\w [dY]\w
[dU]\w [dV]$$ is a $G^J$-invariant volume element on $\BH_{g,h}$,
where
$$[dX]=\w_{\mu\leq\nu}dx_{\mu\nu},\quad [dY]=\w_{\mu\leq\nu}
dy_{\mu\nu},\quad [dU]=\w_{k,l}du_{kl}\quad \text{and} \quad
[dV]=\w_{k,l}dv_{kl}.$$ Using the partial Cayley transform
$\Phi_*$ and Theorem 3.1, we can find a $G_*^J$-invariant
Riemannian metric on the Siegel-Jacobi disk $\BD_{g,h}$ and its
Laplacian explicitly which generalize Formulas (2.12) and (2.13).
For more detail, we refer to \cite{YJH5}.

\newcommand\Imm{\text{Im}\,}
\newcommand\MCM{\mathcal M}

\end{section}

%
%
\begin{section}{{\bf The Canonical Automorphic Factors}}
\setcounter{equation}{0}

The isotropy subgroup $K_*^J$ at $(0,0)$ under the action (3.5) is
\begin{equation}
K^J_*=\left\{\left( \begin{pmatrix} P  &0
\\ 0 &\OP
\end{pmatrix},\,\left(0,0\,;\kappa\right)\right)\,\Big|\ P\in
U(g),\ \kappa\in \BR^{(h,h)}\,\right\}.
\end{equation}
The complexification of $K_*^J$ is given by
\begin{equation}
K^J_{*,\BC}=\left\{\left( \begin{pmatrix} P  &0
\\ 0 &  {}^tP^{-1}
\end{pmatrix},\,\left(0,0\,;\zeta\right)\right)\,\Big|\ P\in
GL(g,\BC),\ \zeta\in \BC^{(h,h)}\,\right\}.
\end{equation}

\noindent By a complicated computation, we can show that if
\begin{equation}
g_*=\left(\begin{pmatrix} P & Q\\
{\overline Q} & {\overline P}
\end{pmatrix},\left( \la, \mu;\,\kappa\right)\right)
\end{equation}
is an element of $G_*^J,$ then the $K_{*,\BC}^J$-component of
$$g_*\cdot\left( \begin{pmatrix} I_g & W\\ 0 & I_g
\end{pmatrix}, (0,\eta;0)\right)$$
is given by
\begin{equation}
\left( \begin{pmatrix} P- (PW+Q)(\OQ W+\OP)^{-1}\OQ &0
\\ 0 &\OQ W+\OP
\end{pmatrix},\,\left(0,0\,;\kappa_*\right)\right),
\end{equation}
where
\begin{eqnarray*}
\ka_*&=&\ka+\la\,\,{}^t\eta+(\eta+\la W+\mu)\,\,{}^t\la-(\eta+\la
W+\mu)\,\,{}^t\OQ\,\,{}^t(\OQ W+\OP)^{-1}\,\,{}^t(\eta+\la
W+\mu)\\
&=&\ka+\la\,\,{}^t\eta+(\eta+\la W+\mu)\,\,{}^t\la-(\eta+\la
W+\mu)(\OQ W+\OP)^{-1}\OQ\,\,{}^t(\eta+\la W+\mu).
\end{eqnarray*}
Here we used the fact that $(\OQ W+\OP)^{-1}\OQ$ is symmetric.
\vskip 0.2cm For $g_*\in G_*^J$ given by (4.3) with $g_0=\begin{pmatrix} P & Q\\
{\overline Q} & {\overline P}
\end{pmatrix}\in G_*$ and
$(W,\eta)\in\BD_{g,h}$, we write
\begin{equation}
J(g_*,(W,\eta))=a(g_*,(W,\eta))\,b(g_0,W),
\end{equation}
where
\begin{equation*}
a(g_*,(W,\eta))=(I_{2g},(0,0;\ka_*)),\quad \textrm{where}\ \ka_* \
\textrm{is\ given in }\ (4.4)
\end{equation*}
and
\begin{equation*}
b(g_0,W)=\left( \begin{pmatrix} P- (PW+Q)(\OQ W+\OP)^{-1}\OQ &0
\\ 0 &\OQ W+\OP
\end{pmatrix},\,\left(0,0\,;0\right)\right).
\end{equation*}
\begin{lemma}
Let
$$\rho:GL(g,\BC)\lrt GL(V_{\rho})$$
be a holomorphic representation of $GL(g,\BC)$ on a finite
dimensional complex vector space $V_{\rho}$ and
$\chi:\BC^{(h,h)}\lrt\BC^{\times}$ be a character of the additive
group $\BC^{(h,h)}.$ Then the mapping
$$J_{\chi,\rho}:G_*^J\times \BD_{g,h}\lrt GL(V_{\rho})$$
defined by
$$J_{\chi,\rho}(g_*,(W,\eta))=\chi(a(g_*,(W,\eta)))\,\rho(b(g_0,W))$$
is an automorphic factor of $G^J_*$ with respect to $\chi$ and
$\rho$.
\end{lemma}
\begin{proof}
We observe that $a(g_*,(W,\eta))$ is a summand of automorphic,
i.e.,
$$a(g_1g_2,(W,\eta))=a(g_1, g_2\cdot(W,\eta))+a(g_2,(W,\eta)),$$
where $g_1,g_2\in G_*^J$ and $(W,\eta)\in\BD_{g,h}.$ Together with
this fact, the proof follows from the fact that the mapping
$$J_{\rho}:G_*\times \BD_g\lrt GL(V_{\rho})$$
defined by
$$J_{\rho}(g_0,W):=\rho(b(g_0,W)),\quad g_0\in G_*,\ W\in\BD_g$$
is an automorphic factor of $G_*.$
\end{proof}

\def\CM{\mathcal M}
\vskip 0.2cm \noindent {\bf Example 4.1.} Let ${\mathcal M}$ be a
symmetric half-integral semi-positive definite matrix of degree
$h$. Then the character
$$\chi_{\CM}:\BC^{(h,h)}\lrt \BC^{\times}$$
defined by
$$\chi_{\mathcal M}(c)=e^{-2\pi i \s(\CM c)},\quad
c\in\BC^{(h,h)}$$ provides the automorphic  factor
$$J_{\CM,\rho}:G_*^J\times \BD_{g,h}\lrt GL(V_{\rho})$$
defined by
$$J_{\CM,\rho}(g_*,(W,\eta))=e^{-2\pi i\s(\CM \ka_*)}\rho(\OQ
W+\OP),$$ where $g_*$ is an element in $G_*^J$ given by (4.3) and
$\ka_*$ is given in (4.4). Using $J_{\CM,\rho}$, we can define the
notion of Jacobi forms on $\BD_{g,h}$ of index $\CM$ with respect
to the Siegel modular group $T^{-1}Sp(g,\BZ)\,T$\,(cf.\,\,[7-9]).

\vskip 0.2cm \noindent $ \textbf{Remark 4.2.}$ The
$P_*^-$-component of
$$g_*\cdot\left( \begin{pmatrix} I_g & W\\ 0 & I_g
\end{pmatrix}, (0,\eta;0)\right)$$
is given by
\begin{equation}
\left( \begin{pmatrix} I_g & 0
\\ (\OQ W+\OP)^{-1}\OQ & I_g
\end{pmatrix},\,\Big(\la-(\eta+\la W+\mu)(\OQ
W+\OP)^{-1}\OQ,\,0\,;0\Big)\right),
\end{equation}
where
\begin{equation*}
P_*^-=\left\{\,\left( \begin{pmatrix} I_g & 0\\ W & I_g
\end{pmatrix}, (\xi,0\,;0)\right)\,\Big|\ W=\,{}^tW\in \BC^{(g,g)},\
\xi\in\BC^{(h,g)}\ \right\}.
\end{equation*}

\end{section}

\end{document}